\documentclass[10pt,reqno]{amsart}
\usepackage{amssymb}
\usepackage[all]{xy}
\newtheorem{thm}{Theorem}[section]
\newtheorem{lem}[thm]{Lemma}
\newtheorem{prop}[thm]{Proposition}

\theoremstyle{definition}
\newtheorem{dfn}[thm]{Definition}

\newtheorem{ex}[thm]{Example}
\theoremstyle{remark}
\newtheorem{remark}[thm]{Remark}

\begin{document}

\vskip 1pc

\title{Topological entropy and AF subalgebras of graph
$C^*$-algebras}

\author[Ja A Jeong]{Ja A Jeong$^{\dagger}$}
\thanks{$^{\dagger}$Partially supported by KOSEF R14-2003-006-01000-0}
\address{
Mathematical Sciences Division\\
BK21, Seoul National University\\
Seoul, 151--742\\
Korea, email: jajeong\-@\-math.\-snu.\-ac.\-kr }
\author[Gi Hyun Park]{Gi Hyun Park$^{\ddagger}$}
\thanks{$^{\ddagger}$Partially supported by KOSEF  R01-2001-000-00001-0}
\address{
Department of Mathematics\\
Hanshin University\\
Osan, 447--791\\
Korea, email: ghpark\-@\-hanshin.\-ac.\-kr }
\date{}
\subjclass[2000]{46L05, 46L55}

\begin{abstract}   Let $\mathcal A_E$ be the canonical AF subalgebra of
a graph $C^*$-algebra $C^*(E)$   associated with a locally finite
directed graph $E$. For Brown-Voiculescu's topological entropy
$ht(\Phi_E)$ of the canonical completely positive map $\Phi_E$ on
$C^*(E)$,  $ht(\Phi_E)=ht(\Phi_E|_{\mathcal A_E})=h_l(E)=h_b(E)$
is known to hold for a finite graph $E$, where $h_l(E)$ is the
loop entropy of Gurevic and $h_b(E)$ is the block entropy of
Salama. For an irreducible infinite graph $E$, the inequality
$h_l(E)\leq ht(\Phi_E|_{\mathcal A_E})$ has been known recently.
It is shown in this paper that
$$ht(\Phi_E|_{\mathcal A_E})\leq \max\{ h_b(E), h_b(\,{}^{t}E)\},$$
where  ${}^tE$ is  the graph $E$ with the direction  of  the edges
reversed. Some Irreducible infinite graphs  $E_p(p>1)$ with
$ht(\Phi_E|_{\mathcal A_{E_p}})=\log p$ are also examined.
\end{abstract}

\maketitle

\setcounter{equation}{0}

\section{Introduction}

 Voiculescu \cite{Vo} introduced a notion of topological
entropy $ht(\alpha)$ for an automorphism $\alpha$ of a nuclear
unital $C^*$-algebra $A$ to measure the growth of $\alpha^n$ as
$n\to \infty$ using the fact that a nuclear $C^*$-algebra has the
completely positive approximation property. The definition extends
very well to automorphisms of exact $C^*$-algebras (as done by
Brown in \cite{Br}) due to the deep result by Kirchberg \cite{Kr}
that exact $C^*$-algebras are nuclearly embeddable. But without
effort one can define $ht(\Phi)$ even for a completely positive
(cp) map on an exact $C^*$-algebra as described in \cite{BG}.
Since a $C^*$-subalgebra of an exact $C^*$-algebra is always
exact, if $\Phi:A\to A$ is a cp map  on an exact $C^*$-algebra $A$
and $B$ is a $\Phi$-invariant $C^*$-subalgebra of $A$ then
$ht(\Phi|_B)$ can be defined and
 the monotonicity $ht(\Phi|_B)\leq
ht(\Phi)$  holds.

  The topological entropy has been computed in several cases, for
  example, the equality $ht(\alpha *\beta)=\max\{ ht(\alpha),
  ht(\beta)\}$ for the reduced free product automorphism
  $\alpha *\beta$ was proved in
  \cite{BDS}, when the free product is with amalgamation over a
  finite dimensional $C^*$-algebra.  Also Dykema  \cite{Dy}  showed that
  $ht(\alpha)=0$ for certain classes of automorphisms $\alpha$ of
  reduced amalgamated free products of $C^*$-algebras, which turns
  out to
  extend St\o rmer's result \cite{St} that the Connes-St\o rmer entropy of
  the free shift automorphism of the II${}_1$-factor $L(F_\infty)$ is
  zero.

In this paper we are concerned  with  the topological entropy of
the shift type cp maps on $C^*$-algebras arising from directed
graphs.
 A typical one  is the canonical cp map $\Phi_A:\mathcal O_A\to \mathcal O_A$
of the Cuntz Krieger algebra   $\mathcal O_A$ given by
$$\Phi_A(x)=\sum_{i=1}^n s_i x
 s_i^*,$$ where $s_1, \dots, s_n$ are the partial isometries that
 generate  $\mathcal O_A$. The reason we call $\Phi_A$ shift type  is
 that  $\mathcal O_A$ contains a $\Phi_A$-invariant commutative
 $C^*$-subalgebra $\mathcal D_A$
which is isomorphic to $C(X_A)$ in such a way
 that the restriction $\Phi_A|_{\mathcal D_A}$ corresponds to the shift map
 $\sigma_{X_A}$ on the (compact) shift space $X_A$ associated
 with the incidence matrix $A$.
 The topological  entropy of $\Phi_A$ is then computed
 (see \cite{Cd}, \cite{BG}, \cite{PWY})
  as $ht(\Phi_A)=\log r(A)$ ($r(A)$ is the spectral radius of  $A$).
  But  $\log r(A)=h_{top}(X_A)$ is a well known
   fact, so that one can deduce by  \cite{De} that
   $ht(\Phi_A)=ht(\Phi_A|_{\mathcal D_A}).$
On the other hand, $\mathcal O_A$ also contains another important
$\Phi_A$-invariant  $C^*$-subalgebra $\mathcal A_A$ which is an AF
algebra with $\mathcal D_A\subset \mathcal A_A$. Thus by
monotonicity of entropy, we have $ht(\Phi_A)=ht(\Phi_A|_{\mathcal
A_A})=ht(\Phi_A|_{\mathcal D_A}).$

The Cuntz-Krieger algebras $\mathcal O_A$ are now well understood
as graph $C^*$-algebras $C^*(E)=C^*(s_e,p_v)$ associated with
finite directed graphs $E$ and  the cp
 map $\Phi_A$ of $\mathcal O_A$  is interpreted as the  map
$\Phi_E:C^*(E)\to C^*(E)$ given by $\Phi_E(x)=\sum_{e\in E^1} s_e
x s_e^*.$ Hence if $E$ is a finite graph (possibly with sinks)
which contains an infinite path, it follows that
$ht(\Phi_E)=ht(\Phi_E|_{\mathcal A_E})=ht(\Phi_E|_{\mathcal
D_E})=\log r(A_E)$, where $\mathcal A_E$ is the AF subalgebra of
$C^*(E)$ corresponding to $\mathcal A_A$ in $\mathcal O_A$ and
$A_E$ is the edge matrix of $E$ (see \cite{JP}).

If $E$ is infinite but locally finite then the map $\Phi_E$ is
known to be a contractive cp map, and furthermore if $E$ is
 irreducible and $\mathcal A_E$ is the canonical AF subalgebra of
 $C^*(E)$, the inequality
 $h_l(E)\leq ht(\Phi_E|_{\mathcal A_E})$ is known to hold \cite{JP}.
   The purpose of
the present paper is then to give an upper bound for
$ht(\Phi_E|_{\mathcal A_E})$ and we actually prove the following
(see Theorem 3.9)
$$ht(\Phi_E|_{\mathcal A_E})\leq \max \{h_b(E),h_b({}^t E)\}.$$
 In particular, for an irreducible infinite
 graph $E_p$ constructed in  \cite{Sa} so that
 $h_l(E_p)=h_b(E_p)=p>1$, we have
   $ht(\Phi_{E_p}|_{\mathcal A_{E_p}})=\log p$.

We believe that the result would be helpful to compute the entropy
$ht(\Phi_E)$ of  $\Phi_E$ on the whole graph $C^*$-algebra
$C^*(E)$.

\vskip 1pc
\section{Preliminaries}

\subsection{Graph $C^*$-algebras}
A (directed) graph is a quadruple $E=(E^0,E^1,r,s)$ of  the vertex
set $E^0$, the edge set $E^1$, and  the range, source maps
$r,s:E^1\to E^0$. A family $\{p_v,s_e\mid v\in E^0, e\in E^1 \}$
of mutually orthogonal projections $p_v$ and partial isometries
$s_e$  is called a {\it Cuntz-Krieger $E$-family} if the following
relations
hold; \begin{align*} &s_e^*s_e=p_v, \ \ s_es_e^*\leq p_{s(e)},\\
&p_v=\sum_{s(e)=v} s_es_e^*, \ \ \text{ if }
0<|s^{-1}(v)|<\infty.\end{align*} The {\it graph $C^*$-algebra}
$C^*(E)$ is then defined to be a $C^*$-algebra generated by a
universal Cuntz-Krieger $E$-family (see \cite{KPR}, \cite{BPRS}).
We call $E$  {\it locally finite} if each vertex receives and
emits only finitely many edges. Throughout this paper we consider
only locally finite graphs, and adopt the  notations  in
\cite{KPR}.
 If a finite path $\alpha\in E^*$ of length $|\alpha|>0$ is a
 return path, that is,
 $s(\alpha)=r(\alpha)$, then $\alpha$ is called a {\it loop} at $v=s(\alpha)$.
  A graph $E$ is said to be {\it irreducible} if for any two vertices $v,w$
there is a finite path $\alpha\in E^*$ with $s(\alpha)=v$ and
$r(\alpha)=w$. It is known that if $E$ is irreducible and every
loop has an exit then $C^*(E)$ is simple.

\subsection{Topological entropy of cp maps}

    Let $A$ be a $C^{*}$-algebra, $\pi:A \to B(H)$  a faithful
    $*$-representation, and  $Pf(A)$ be the set of all finite
    subsets of $A$. For $\omega\in Pf(A)$ and $\delta>0$, put
 \begin{align*}
 CPA(\pi,A)& = \{ (\phi, \psi, B) \mid \phi : A\to B, \psi : B \to B(H)
    \text{ are contractive cp maps }\\
  &\ \hskip 2cm  \text{  and } \dim B < \infty \},\\
    rcp(\pi,\omega,\delta)& =\inf \{ rank(B) \mid (\phi,\psi, B) \in
    CPA(\pi, A),\ \| \psi \circ \phi (x) - \pi (x) \| < \delta,\\
    & \ \hskip 2cm
    \text{ for all } x \in \omega \},\end{align*}
 where $rank(B)$ denotes the dimension
    of a maximal abelian subalgebra of $B$.

Since the cp $\delta$-rank $rcp(\pi,\omega,\delta)$ is independent
of the choice of $\pi$ (see \cite{Br}, \cite{BG}) and graph
$C^*$-algebras $C^*(E)$ are nuclear we may write $rcp(
\omega,\delta)$ for $rcp(\pi, \omega,\delta)$ assuming that
$C^*(E)\subset B(H)$  for a Hilbert space $H$.\vskip 1pc

\begin{dfn} \rm (\cite{Br}, \cite{BG})
    Let $A\subset B(H)$ be a $C^*$-algebra and  $\Phi:A\to A $
be  a cp map. Put
 \begin{align*} ht(\Phi,\omega,\delta) &=
\limsup_{n\to \infty} \frac{1}{n} \log\big( rcp(\,\cup_{i=0}^{n-1}
\Phi^{i}(\omega), \delta)\big)\\
ht(\Phi,\omega)&=\sup_{\delta>0}
ht(\Phi,\omega,\delta).\end{align*}
     Then $ht(\Phi):=\sup_{\omega\in Pf(A)}
ht(\Phi,\omega)$ is called the {\it topological entropy} of
 $\Phi$.
\end{dfn}

\vskip .5pc
\begin{remark} \rm We refer the reader to \cite{BG} and \cite{Br}
for the following useful properties  and  their proofs. Let $A$ be
an exact $C^*$-algebra and $\Phi:A\to A$ be a cp map.

\renewcommand{\theenumi}{\arabic{enumi}}
\renewcommand{\labelenumi}{{\theenumi}.}
\setlength{\itemsep}{10mm} \begin{enumerate}\item[(a)] If
 $\theta:A\to B$ is a
$C^*$-isomorphism then  $ht(\Phi)=ht(\theta\Phi\theta^{-1}).$

\item[(b)] Let $\tilde A$ be the unital $C^*$-algebra obtained by
adjoining a unit and $\tilde \Phi : \tilde A \to \tilde A$ be the
extension of $\Phi$. Then $ht(\tilde \Phi) =ht(\Phi).$

\item[(c)] If $A_0 \subset A$ is a $\Phi$-invariant
$C^*$-subalgebra of $A$, $ht( \Phi|_{A_0}) \leq ht(\Phi).$

\end{enumerate}
\end{remark}

 We will use the following Arveson's extension
theorem several times.

\vskip .5pc

\noindent {\bf Arveson Extension Theorem} (see \cite{Br}) Let $A$
be a unital $C^*$-algebra, $S\subset A$  a unital subspace with
$S=S^*$, and $\phi:S\to B$ be a contractive cp map where $B=B(H)$
or $dim(B)<\infty$. Then $\phi$ extends to a cp map
$\bar{\phi}:A\to B$. If $S$ is a $C^*$-subalgebra of $A$ then we
obtain a unital cp extension of $\phi$ even when $S$ does not
contain the unit of $A$.

\vskip 1pc

If $E$ is a locally finite graph, the map  $\Phi_E:C^*(E)\to
C^*(E)$, defined by
$$ \Phi_E(x)= \sum_{e\in E^1} s_e xs_e^*,$$ is well defined,
contractive, and completely positive \cite{JP}.
 For a finite
graph $E$, the topological entropy $ht(\Phi_E)$ has been obtained
as follows (see \cite{BG}, \cite{Cd}, \cite{PWY}, or  \cite{JP}).

\vskip 1pc
\begin{thm}
Let $E$ be a finite  graph possibly with sinks and $A_E$ be the
edge matrix of $E$. If  $E$ contains an infinite path  then
$$ht(\Phi_E)=\log r(A_E),$$ where $r(A_E)$ is the
spectral radius of $A_E$.
\end{thm}
\vskip .5pc

 By $h_{top}(X)$ we denote the topological entropy of a
compact space $(X,T)$ together with a continuous map $T:X\to X$
(for definition, see \cite[Definition 4.1.1]{LM} or
\cite[p.23]{Kt} ). Let $E$ be a locally finite infinite graph and
 $X_E$ the locally compact shift space of (one-sided) infinite paths
 with the one point compactification $\bar{X}_{E}$.
Then the first identity in the following theorem is shown for the
doubly infinite path space of $E$ by Gurevic \cite{Gu}. See
Definition 3.1 for $h_b(E)$.

\begin{thm}[\expandafter{\cite[Theorem 4.4]{JP}}]
Let $E$ be a locally finite irreducible infinite graph. Then
$$ h_{top}(\bar{X}_{E})=\sup_{E'}h_b(E')\,\leq \, ht(\Phi_E),$$
where the supremum  is taken over all the finite (irreducible)
subgraphs of $E$.
\end{thm}

\vskip 1pc

\section{Main results}
 Throughout this section $E$ will denote a locally
finite infinite graph unless stated otherwise.
 For a path $\alpha\in E^*$, let $\alpha^0$ be the set of
vertices lying on $\alpha=\alpha_1\cdots\alpha_n$, that is,
$\alpha^0=\{s(\alpha_1), r(\alpha_1), \dots , r(\alpha_n)\}$. For
a fixed vertex $v$  we consider the following subsets of finite
paths $E^n$ of length $n$.

\renewcommand{\theenumi}{\arabic{enumi}}
\renewcommand{\labelenumi}{{\theenumi}.}
\setlength{\itemsep}{10mm}
\begin{enumerate}
 \item [(i)]$E^n(v)=\{\alpha\in E^n\mid v\in \alpha^0\}$,
 \item[(ii)] $E_s^n(v)=\{\alpha\in E^n\mid s(\alpha)=v\}$,
 \item[(iii)] $E_{s}^n(v^\star)=\{\alpha \in E_s^n(v)\mid
r(\alpha_i)\neq v, \  1\leq i\leq n\}$,
 \item[(iv)] $E_l^n(v)=\{\alpha\in E^n\mid \alpha \text{ is a loop at }
v\}$.\end{enumerate} Similarly we can think of $E_r^n(v)$ and
$E_{r}^n(v^\star)$.

\renewcommand{\theenumi}{\arabic{enumi}}
\renewcommand{\labelenumi}{{\theenumi}.}
\setlength{\itemsep}{10mm}

\vskip 1pc \begin{dfn} Let $E$ be an irreducible graph and let
$v\in E^0$.
\begin{enumerate}
 \item[(a)]
 $h_l (E):=\limsup_n \frac{1}{n}\log |E_l^n(v)|$
 is called the {\it loop
entropy} of $E$.

 \vskip .5pc\noindent \item[(b)]  $h_{b} (E):=\limsup_n \frac{1}{n}
\log |E_s^n(v)|$ is called the {\it block entropy} of $E$.
\end{enumerate}\end{dfn}

\noindent Note  that both entropies $h_l(E)$ and  $h_b(E)$
 are independent of the choice of a
vertex $v$ \cite{Sa}.
    If ${}^tE$  denotes  the graph $E$ with the
direction of all edges  reversed, then clearly $h_l (E)=h_l
({}^tE)$ while $h_b(E)\neq h_b({}^t E)$ in general as we will see
in Example 3.3.

We will use the following notation for the infinite series with
coefficients from (i)-(iv) above.
\begin{enumerate}
\item[(i)$'$] $E(v,z):=\sum |E^n(v)|z^n,$
 \item[(ii)$'$] $E_s(v,z):=\sum |E_s^n(v)| z^n,$
\item[(iii)$'$] $E_{s}(v^\star,z):=\sum |E_{s}^n(v^\star)| z^n,$
\item[(iv)$'$] $E_l(v,z):=\sum |E_l^n(v)|z^n.$ \end{enumerate}

 We denote the
radius of convergence of the series $E_s(v^\star,z)$ by
$R_{E_s^\star}$. Thus
$$R_{E_s^\star}^{-1}=\limsup_{n\to \infty} \big|
E_s^n(v^\star)\big|^{1/n}.$$ Similarly $R_{E_r^\star}$ denotes the
radius of convergence of $E_r(v^\star,z):=\sum |E_{r}^n(v^\star)|
z^n$.

 \vskip 1pc
\begin{prop}[\cite{Sa}] If $E$ is an irreducible graph, then
 $$h_b(E)=\max\{ \log \big( R_{E^\star_{s}}^{-1}\big), \, h_l(E)\}.$$
 \end{prop}

\vskip 1pc \noindent  Note that if $E$ is irreducible then
 $h_b({}^tE)=\limsup \frac{1}{n} \log |E_r^n(v)|$ and so
  from the above proposition we have
$$h_b({}^tE)=\max\{\log \big( R_{E^\star_{r}}^{-1}\big), \, h_l(E)\}.$$

\vskip 1pc
 The following example shows that $h_b(E)\neq h_b({}^tE)$ in
 general.
\vskip 1pc
\begin{ex} For each pair of positive real numbers
$1<p\leq q$, \rm Salama \cite{Sa} constructed an irreducible
infinite graph $E_{p,q}$ with $$h_l(E_{p,q})=\log p\ \ \text{and}\
\ h_b(E_{p,q})=\log q.$$  For example,  the following graph
$E:=E_{2,8}$ satisfies $h_l(E)=\log 2$ and $h_b(E)=\log 8$.
 There are 8 edges from the vertex $n$ to the vertex $n+1$ for each $n\geq
 0$.

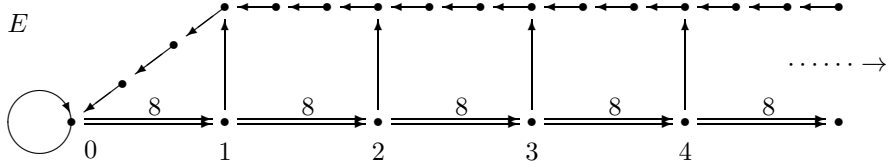
\setlength{\unitlength}{1.7cm}

\hspace*{6cm}
\begin{picture}(8,1.5)

\put(-3.25,0){\circle{0.5}}\put(-2.9,-0.28){$0$}
\put(-3,0){\circle*{0.07}}
 \put(-2.6,0.3){\circle*{0.07}} \put(-2.2,0.6){\circle*{0.07}}
\put(-1.8,0){\circle*{0.07}}\put(-1.85,-0.3){1}
\put(-0.6,0){\circle*{0.07}}\put(-0.65,-0.3){$2$}
\put(0.6,0){\circle*{0.07}}\put(0.55,-0.3){$3$}
\put(1.8,0){\circle*{0.07}}\put(1.75,-0.3){$4$}

\put(3,0){\circle*{0.07}}

 \put(-1.8,0.9){\circle*{0.07}}
\put(-1.4,0.9){\circle*{0.07}}\put(-1,0.9){\circle*{0.07}}

\put(-0.6,0.9){\circle*{0.07}}\put(-0.2,0.9){\circle*{0.07}}\put(0.2,0.9){\circle*{0.07}}
\put(0.6,0.9){\circle*{0.07}}
\put(1,0.9){\circle*{0.07}}\put(1.4,0.9){\circle*{0.07}}\put(1.8,0.9){\circle*{0.07}}
\put(2.2,0.9){\circle*{0.07}}\put(2.6,0.9){\circle*{0.07}}\put(3,0.9){\circle*{0.07}}

\put(-1.4,0.9){\vector(-1,0){0.3}}
\put(-1,0.9){\vector(-1,0){0.3}}
\put(-0.6,0.9){\vector(-1,0){0.3}}
\put(-0.2,0.9){\vector(-1,0){0.3}}
\put(0.2,0.9){\vector(-1,0){0.3}}
\put(0.6,0.9){\vector(-1,0){0.3}}

\put(1.0,0.9){\vector(-1,0){0.3}}

\put(1.4,0.9){\vector(-1,0){0.3}}

\put(1.8,0.9){\vector(-1,0){0.3}}

\put(2.2,0.9){\vector(-1,0){0.3}}

\put(2.6,0.9){\vector(-1,0){0.3}}

\put(3,0.9){\vector(-1,0){0.3}}

\put(-3.01,0.07){\vector(1,-2){0}}

\put(-2.9,0.02){\vector(1,0){1.0}}\put(-2.9,-0.02){\vector(1,0){1.0}}

\put(-1.7,0.02){\vector(1,0){1.0}}\put(-1.7,-0.02){\vector(1,0){1.0}}

\put(-0.5,0.02){\vector(1,0){1.0}}\put(-0.5,-0.02){\vector(1,0){1.0}}

\put(0.7,0.02){\vector(1,0){1.0}}\put(0.7,-0.02){\vector(1,0){1.0}}

\put(1.9,0.02){\vector(1,0){1.0}}\put(1.9,-0.02){\vector(1,0){1.0}}

\put(-1.8,0.1){\vector(0,1){0.7}}

\put(-0.6,0.1){\vector(0,1){0.7}}

\put(0.6,0.1){\vector(0,1){0.7}}

\put(1.8,0.1){\vector(0,1){0.7}}

\put(-2.6,0.3){\vector(-4,-3){0.3}}

\put(-2.2,0.6){\vector(-4,-3){0.3}}

\put(-1.8,0.9){\vector(-4,-3){0.3}}

\put(2.6,0.4){$\cdots\cdots\rightarrow$}

\put(-2.4,0.05){8}

\put(-1.2,0.05){8}\put(0,0.05){8}\put(1.2,0.05){8}\put(2.4,0.05){8}

\put(-3.5,0.7){$E$}
\end{picture}

\vskip 1cm \noindent Now we show that $$\log \big(
R_{E^\star_{r}}^{-1}\big) \leq h_l(E),$$ which then implies
$h_b({}^tE)=h_l(E)$ by the above proposition (hence
$h_b({}^tE)\neq h_b(E)$). For a fixed vertex $0$ we have
\begin{align*}R_{E^\star_r}^{-1}&=\limsup_{n\to \infty}|E_{r}^{n}(0^\star)|^{1/n}\\
&=\limsup_{n\to \infty}\big| \{\alpha\in E_r^n(0)\mid
s(\alpha_i)\neq 0, \text{ for } 1\leq i\leq n
\}\big|^{1/n}.\end{align*}
 With  $ n_{k+1}=4k+1\, (k\geq 0),$
$$|E_{r}^{n_{k+1}-1}(0^\star)|=|E_{r}^{4k}(0^\star)|
 =1+8^{k-1}+8^{k-4}+8^{k-7}+\cdots,$$
 and a computation gives
  $$\limsup_{k\to \infty} |E_{r}^{4k}(0^\star)|^{\frac{1}{4k}}= 8^{1/4}.$$
 But it is not hard to see that $$\limsup_{n\to \infty} |E_{r}^{n}(0^\star)|^{\frac{1}{n}}=
  \limsup_{k\to \infty} |E_{r}^{4k}(0^\star)|^{\frac{1}{4k}}\, ,$$
hence $\log (R_{f^\star_r}^{-1})=\log 8^{1/4}<\log2=h_l(E).$
\end{ex}

\vskip 1pc

\begin{lem} If $E$ is an irreducible graph then the value
$$ \limsup_{n\to \infty} \frac{1}{n} \log |E^n(v)|
$$ is independent
of the choice of a vertex $v$.
\end{lem}

\proof
 Let
$v,w$ be two vertices of $E$. Then there exist two paths $\mu\in
E^{k}$,
  $\nu\in E^{m}$ with $s(\mu)=r(\nu)=v$, $s(\nu)=r(\mu)=w$
  because $E$ is irreducible. We assume that $\mu$ and $\nu$ have
  the smallest length, respectively.
  If $\alpha=\alpha_1\alpha_2\cdots \alpha_n\in E^n(v)$ then with $i_0=\min \{i\mid
  s(\alpha_i)=v\}$
  write $\alpha=\alpha'\alpha''$, where
  $\alpha'=\alpha_1\cdots\alpha_{i_0-1}$ and
  $\alpha''=\alpha_{i_0}\cdots\alpha_n$ (if $i_0=1$,
  $\alpha=\alpha''$).
  Then the map $$E^n(v)\to E^{n+k+m}(w),\ \
  \alpha=\alpha'\alpha''\mapsto \alpha'\mu\nu\alpha''$$ is
  injective, hence $|E^n(v)|\leq |E^{n+k+m}(w)|$ for each
  $n$. Therefore \begin{align*}
  \limsup_{n\to \infty} \frac{1}{n} \log |E^n(v)|&\leq
  \limsup_{n\to \infty} \frac{1}{n} \log |E^{n+k+m}(w)|\\&\leq
  \limsup_{n\to \infty} \frac{1}{n} \log |E^n(w)|.\end{align*}
\endproof

\medbreak \vskip 1pc

\begin{prop} Let $E$ be an irreducible graph and $v_0\in E^0$.

\noindent {\rm (a)} If $E$ is  finite, then $$\limsup_{n\to
\infty} \frac{1}{n} \log |E_l^n(v_0)|=\limsup_{n\to \infty}
\frac{1}{n} \log |E^n|.$$

In particular, $h_l(E)=h_{b}(E)=h_{b}({}^tE).$

\noindent {\rm (b)} If $E$ is  infinite,  then
$$ \limsup_{n\to \infty} \frac{1}{n} \log |E^n(v_0)|
\ = \ \max \{h_b(E), \, h_b({}^tE)\}.$$
\end{prop}

\proof (a)  Let $E^0=\{v_0, v_1, \dots, v_{k-1}\}$. Since $E$ is
irreducible there exist finite paths $\{\mu_i,\, \nu_i\mid 0\leq
i\leq k-1\}$
 such that $ s(\mu_i)=r(\nu_i)=v_0$, $r(\mu_i)=v_i=s(\nu_i).$
 Suppose $|\mu_i|=m_i$, $|\nu_j|=l_j.$
     If $\alpha\in E^n$ is a path with $s(\alpha)=v_i$,
     $r(\alpha)=v_j$ then $\mu_i\alpha \nu_j\in E_l^{n+m_i+l_j}(v_0)$
     is a loop at $v_0$.
 The map $\alpha\mapsto \mu_i\alpha \nu_j$ is not necessarily injective, but
 there exist at most $k_0$ paths in $E^n$ that have the same image in
 $E_l^{n+m_i+l_j}(v_0)$ under the map, where $k_0=\max_{i,j}\{m_i+l_j\}$.
 Hence we  have
$$ |E^n|\leq k_0\cdot \big|\cup_{0\leq i,j\leq k-1}
 E_l^{n+m_i+l_j}(v_0)\big|
 \leq k_0 k^2 \max_{i,j}|E_l^{n+m_i+l_j}(v_0)|.$$
On the other hand, for each $n$, there exists a $k_n\in \{0,\dots
,k_0\}$ such that
$$|E_l^{n+k_n}(v_0)|=\max_{i,j}|E_l^{n+m_i+l_j}(v_0)|.$$
Then $|E^n|\leq k_0k^2|E_l^{n+k_n}(v_0)|$ and it follows that
$$\limsup_{n\to \infty}\frac{1}{n}\log |E^n| \leq \limsup_{n\to
\infty}\frac{1}{n}\log |E_l^n(v_0)|.$$

\noindent (b) Note first that
\begin{align*}
|E^n(v)| &= \big| \cup_{k=0}^n\{\alpha\beta\mid \alpha\in
E_{r}^k(v^\star),\, \beta\in E_s^{n-k}(v_0)\}\big|\\&=
\sum_{k=0}^n \big|E_{r}^k(v^\star)\big|\, \big| E_s^{n-k}(v)
\big|=\sum_{k=0}^n\big| ({}^tE)_{s}^k(v^\star)\big|\, \big|
E_s^{n-k}(v) \big|.\end{align*} Then
 \begin{align*}E(v,z)
 &=\sum_n \Big(\sum_{k=0}^n |({}^tE)_{s}^k(v^\star)|\, | E_s^{n-k}(v)
 |\Big)z^n\\
 &=\Big(\sum_n |({}^tE)_{s}^n(v^\star)|z^n \Big)\, \Big( \sum_n | E_s^{n}(v)|z^n\Big)\\
 &=({}^tE)_{s}(v^\star,z)\cdot E_s(v,z),\end{align*}
so that the radius  of convergence $R_E$ of $E(v,z)$ is equal to
  $\min \big\{ R_{({}^tE)^\star_s}, \,
 R_{E_s}\big\}.$
  Thus  $$R_E^{-1}=\max \big\{\,  R_{({}^tE)^\star_{s}}^{-1}, \,
 R_{E_s}^{-1}\,\big\}.$$
  But  Proposition 3.2 gives $$\log \big(R_{({}^tE)^\star_s}^{-1}\big)\leq
 h_b({}^tE),$$ and also by definition $\log (R_{E_s}^{-1})=h_b(E).$
 Therefore
 $$\limsup_{n\to \infty} \frac{1}{n} \log |E^n(v)|
 =\log (R_E^{-1})\leq \max \{h_b({}^tE),\, h_b(E)\}.$$
\endproof  \medbreak

\vskip 1pc

Let $E$ be an irreducible infinite graph and let $\mathcal D_E$ be
the commutative $C^*$-subalgebra of $C^*(E)$ generated by the
projections $\{p_\alpha=s_\alpha s_\alpha^*\mid \alpha\in E^*\}$.
Then $\mathcal D_E=\overline{\rm span}\{p_\alpha\mid \alpha\in
E^*\}$ and    the map
$$w: \mathcal D_E\to C_0(X_E),\ w(p_\alpha)=\chi_{[\alpha]},$$  is
a $C^*$-isomorphism \cite{JP}. Here $\chi_{[\alpha]}$ is the
characteristic function  on the cylinder  set $[\alpha]=\{\beta\in
X_E\mid \beta=\alpha\beta'\}$ which is both open and closed.
Furthermore, from the proof of Theorem 2.4 we know that
$h_{top}(\bar{X}_E)=ht(\Phi_E|_{\mathcal D_E}).$
   Put $$\mathcal A_E:=\overline{\rm span}\{\, s_\alpha s_\beta^*\,\mid\,
\alpha,\beta\in E^*, \ |\alpha|=|\beta|\,\}.$$ Then  $\mathcal
A_E$ is a $\Phi_E$-invariant AF $C^*$-subalgebra of $C^*(E)$ with
$\mathcal D_E\subset \mathcal A_E$, hence
$$ht(\Phi_E|_{\mathcal D_E})\leq ht(\Phi_E|_{\mathcal A_E}).$$
  Also, for each $v\in E^0$, $\mathcal
A_E$ contains a  $\Phi_E$-invariant AF subalgebra  $\mathcal
A_E(v)$, $$\mathcal A_E(v):=\overline{\rm span}\{\, s_\alpha
s_\beta^*\,\mid\, \,r(\alpha)=r(\beta)=v,\, |\alpha|=|\beta|\,
\}.$$

\vskip 1pc

\begin{lem} Let $v$ be a vertex of an irreducible  graph $E$ with
at least two vertices
 and let $ n\geq 1$.  Then the elements in the set $$\omega(n,v)=\big\{s_\alpha
s_\beta^*\, \big|\ \, r(\alpha)=r(\beta)=v,\, |\alpha|=|\beta|\leq
n \big\}$$ are linearly independent.
\end{lem}
\proof We prove the assertion by induction on $n$.
  For $n=1$, suppose $$x=\sum_{\substack{{e,f\in E^1}\\
  {r(e)=r(f)=v}}} \lambda_{ef}s_e s_f^*+\lambda_0 p_v=0.$$
  If $e_0$ and $f_0$ are edges with $r(e_0)=r(f_0)=v$ and
 either $s(e_0)\neq v$ or $s(f_0)\neq v$  then $s_{e_0}^*p_v
 s_{f_0}=0$, hence
$$0=s_{e_0}^* x s_{f_0}=\lambda_{e_0 f_0}(s_{e_0}^*s_{e_0}
)(s_{f_0}^*s_{f_0})=\lambda_{e_0 f_0}p_v,$$ thus $\lambda_{e_0
f_0}=0$. Similarly, $\lambda_{e f}=0$ if $e$ and $f$ are loops at
$v$ and $e\neq f$. Then $x$ becomes $$x=\sum_{e\in E_l^1(v)}
\lambda_{ee} s_e s_e^* +\lambda_0 p_v=0.$$
  By  irreducibility of $E$ and the assumption that $|E^0|>1$,
  there exists an edge $f$ with
  $s(f)=v, r(f)\neq v$. Then $s_f s_f^* x=\lambda_0 s_f
  s_f^*=0$, so that $\lambda_0=0$ and we have $x=\sum_{e\in E_l^1(v)}
   \lambda_{ee} s_e s_e^*=0.$ Since the projections
  $\{s_e s_e^*\mid e\in E_l^1(v)\}$  are mutually orthogonal, it
  follows that $\lambda_{ee}=0$ for each $e\in E_l^1(v)$.

  Now suppose that the assertion is
true for $n-1$. If
$$x=\sum_{\substack{{|\alpha|=|\beta|\leq
n}\\{r(\alpha)=r(\beta)=v}}} \lambda_{\alpha\beta} s_\alpha
s_\beta^*=0,\ \lambda_{\alpha\beta}\in \mathbb C,$$ then for an
edge $e\in E^1$ we have $$0=s_e^* x
s_e=\sum_{\substack{{\alpha=e\alpha'}\\{\beta=e\beta'}}}\lambda_{\alpha\beta}
s_e^* s_\alpha s_\beta^* s_e= \sum_{|\alpha'|=|\beta'|\leq n-1}
\lambda_{(e\alpha')(e\beta')} s_{\alpha'}(s_{\beta'})^*.$$
 Note that the elements $s_{\alpha'}(s_{\beta'})^*$ appearing in the sum
 are distinct. Thus
by induction hypothesis, one sees that
$\lambda_{(e\alpha')(e\beta')}=0$. But the edge $e$ was arbitrary,
and so we conclude that the coefficients $\lambda_{\alpha\beta}$
are all zero. \endproof  \medbreak

\vskip 1pc

\begin{prop} {\rm (cf. \cite[Proposition 2.6]{Br})}
Let $\Phi:A\to A$ be a  contractive cp map of an exact
$C^*$-algebra $A$. If $\{\omega_\lambda\}_{\lambda\in \Lambda}$ is
a net (partially ordered by inclusion) of finite subsets in $A$
such that the linear span of $\ \cup_{\lambda,l\in \mathbb
Z^+}\Phi^l(\omega_\lambda)$ is dense in $A$ then
$$ht(\Phi)=\sup_\lambda ht(\Phi, \omega_\lambda ).$$
\end{prop}

\vskip 1pc

\begin{thm}
 Let $E$ be an irreducible infinite graph. Then for each $v\in
 E^0$,
 $$ht(\Phi_E|_{\mathcal A_E(v)})\,\leq\,
 h_b({}^tE).$$
 \end{thm}

\proof  Let  $A_n(v)$ be the $C^*$-subalgebra of $\mathcal A_E(v)$
generated by  $\omega(n,v).$
 Then  from
\[ s_\alpha s_\beta^*\cdot s_\mu s_\nu^*=\left\{
\begin{array}{ll} s_{\alpha\mu'}s_\nu^*, & \text{if }\ \mu=\beta\mu',\\
s_\alpha s_{\nu\beta'}^*, & \text{if}\ \beta=\mu\beta',\\
0, & \text{otherwise,}
 \end{array} \right.   \] we see that $A_n(v)={\rm span}(\omega(n,v))$
 is finite dimensional.

 Since $\{\omega(n,v)\}_n$ is an increasing
sequence of finite subsets in $\mathcal A_E(v)$ such that  the
linear span of $\cup_n \omega(n,v)$ is dense in $\mathcal A_E(v)$,
 by Proposition 3.7 it suffices to show that
$$ ht(\Phi_E, \omega(n,v))\leq  h_b({}^tE), \ \ n\in \mathbb N.$$
    Set
$E_l^*(v):=\cup_{k\geq 0} E_l^k(v)$ and  $r(n):=|\,\cup_{k=0}^n
E_r^k(v)\,|.$
   Fix $n_0\in \mathbb N$, and define a map
$\phi:\omega(n_0,v)\to M_{r(n_0)}$  by
$$ \phi (s_\alpha
s_\beta^*)= \sum_{\substack{{|\alpha\gamma|\leq n_0}\\{\gamma\in
E_l^*(v)}}} e_{(\alpha\gamma)(\beta\gamma)},$$ where
$\{e_{\mu\nu}\}$ is the standard matrix units of of the matrix
algebra $M_{r(n_0)}$.
  Since the elements in $\omega({n_0},v)$ are linearly
independent by Lemma 3.6, one can extend the map to the linear map
$\phi: A_{n_0}(v)\to M_{r(n_0)}$.  Now we show that $\phi$ is in
fact a $*$-isomorphism. To prove that it is a $*$-homomorphism, we
only need to see that
$$\phi((s_\alpha s_\beta^*)(s_\mu s_\nu^*))=
\phi(s_\alpha s_\beta^*)\phi(s_\mu s_\nu^*).$$
   If
$\beta=\mu\beta'$,  then $s_\alpha s_\beta^* s_\mu
s_\nu^*=s_\alpha (s_{\nu\beta'})^*$
   and \begin{align*}\phi(s_\alpha
s_\beta^*)\phi(s_\mu s_\nu^*)&=\sum_{\substack{{|\alpha\gamma|\leq
n_0}\\{\gamma\in E_l^*(v)}}}
e_{(\alpha\gamma)(\mu\beta'\gamma)}\cdot
\sum_{\substack{{|\mu\delta|\leq n_0}\\{\delta\in E_l^*(v)}}}
e_{(\mu\delta)(\nu\delta)}\\
 &= \sum_{\substack{{|\alpha\gamma|\leq n_0}\\{\gamma\in E_l^*(v)}}}
e_{(\alpha\gamma)(\nu\beta'\gamma)}=\phi(s_\alpha
(s_{\nu\beta'})^*)=\phi(s_\alpha s_\beta^* s_\mu s_\nu^*).
\end{align*}
  If $\mu=\beta\mu'$, a similar proof works. Otherwise,
$\phi((s_\alpha s_\beta^*)(s_\mu s_\nu^*))=0= \phi(s_\alpha
s_\beta^*)\phi(s_\mu s_\nu^*).$
    In order to show that $\phi$ is injective, let
  $\phi(\sum_{\alpha,\beta}\lambda_{\alpha\beta}\, s_\alpha
  s_\beta^*)=0$. Then
  $$\sum_{\alpha,\beta}\lambda_{\alpha\beta}\phi(s_\alpha
  s_\beta^*)=\sum_{\alpha,\beta}\lambda_{\alpha\beta}
  \big(\sum_{\substack{{|\alpha\gamma|\leq {n_0}}\\{\gamma\in E_l^*(v)}}}
e_{(\alpha\gamma)(\beta\gamma)}\big)=0.$$
 But  the vectors,
 $\sum_{\substack{{|\alpha\gamma|\leq {n_0}}\\{\gamma\in E_l^*(v)}}}
e_{(\alpha\gamma)(\beta\gamma)}$ $(r(\alpha)=r(\beta)=v,\,
|\alpha|=|\beta|\leq {n_0}),$ are linearly independent in
$M_{r(n_0)}$. In fact, if
$A:=\sum_{\alpha,\beta}\lambda_{\alpha\beta}
  \big(\sum_{\substack{{|\alpha\gamma|\leq {n_0}}\\{\gamma\in E_l^*(v)}}}
e_{(\alpha\gamma)(\beta\gamma)}\big)=0,$ then
$e_{vv}Ae_{vv}=\lambda_{vv}e_{vv}=0$, that is, $\lambda_{vv}=0$,
and for any $\alpha,\beta\in E^1_r(v)$,
$e_{\alpha\alpha}Ae_{\beta\beta}=\lambda_{\alpha\beta}e_{\alpha\beta}=0$,
hence  $\lambda_{\alpha\beta}=0.$
   Repeating the process one has
$\lambda_{\alpha\beta}=0$ for any $\alpha,  \beta\in
\cup_{k=0}^{n_0} E^k_r(v).$
   Therefore
$\sum_{\alpha,\beta}\lambda_{\alpha\beta}\, s_\alpha s_\beta^*=0$,
and  the map $\phi$ is injective. The surjectivity of $\phi$
follows from
 $\dim(A_{n_0}(v))=r(n_0)^2$. We simply write $\phi$ for
   $\phi:A_{n_0+l}(v)\to M_{r(n_0+l)}$ $(l\geq 0)$,
   and $\bar\phi$ for its contractive cp
 extension to $\mathcal A_E(v)$ that exists by   Arveson's extension
 theorem.

 For each $n\in \mathbb N$ and  $0\leq l\leq n-1$, note that
  $$\cup_{l=0}^{n-1}\Phi_E^l(\omega(n_0,v))\,\subseteq \,{\rm
 span}(\omega(n_0+n-1,v)).$$
 Then the element
 $$(\bar\phi,\,\psi:=\phi^{-1}\, ,M_{r(n_0+n-1)})
 \in CPA(id,\mathcal A_E(v))$$ satisfies
 $\psi\circ\bar\phi|_{\omega(n_0+n-1,v)}=id_{\omega(n_0+n-1,v)}.$ Thus for
 each $\delta>0$
$$rcp(id,\omega(n_0+n-1,v),\delta)\leq r(n_0+n-1),$$ and so
\begin{align*}ht(\Phi_E|_{\mathcal A_E(v)}, \omega(n_0,v),\delta)&\leq
\limsup_{n\to \infty}\frac{1}{n}\log (r(n_0+n-1))\\
&=\limsup_{n\to \infty}\frac{1}{n}\log (r(n))\\
&=\limsup_{n\to \infty}\frac{1}{n}\log \big|\,\cup_{k=0}^n E_r^k(v)\,\big|\\
&=h_b({}^tE).\end{align*} For the last equality, note that if
$k\leq n$ then $ \big|\, E_r^k(v)\,\big|\leq \big|\,
E_r^n(v)\,\big|,$
 hence $\big|\,\cup_{k=0}^n E_r^k(v)\,\big|\leq (n+1)\cdot
 \big|\,E_r^n(v)\,\big|.$
\endproof  \medbreak

\vskip 1pc Since  $h_l (E)=\sup_{E'\subset E}h(X_{E'})$ is well
known (see \cite[Proposition 7.2.6]{Kt}), we see from Theorem 2.4
and its proof that $h_l(E)\leq ht(\Phi_E|_{\mathcal A_E})$ holds.
The following theorem gives an upper bound for
$ht(\Phi_E|_{\mathcal A_E})$.

\vskip 1pc
\begin{thm} Let $E$ be an irreducible infinite graph and $\mathcal A_E$ be
the AF subalgebra of  $C^*(E)$ generated by the partial isometries
$\{s_\alpha s_\beta^*\mid \alpha, \beta\in E^*,
|\alpha|=|\beta|\}$. Then
$$ht(\Phi_E|_{\mathcal A_E}) \leq \ \max\{h_b({}^t E), h_b(E)\}.$$
\end{thm}

\proof  Let $E^0=\{v_1, v_2, \cdots\}$. For each $n_0\in \mathbb
N$ and $n_1\in \mathbb Z^+=\{0\}\cup \mathbb N$, put
\begin{align*}
\omega(n_0,n_1)&:=\Big\{s_\alpha s_\beta^*\,\big|\,
\alpha,\beta\in
E^{n_1}, r(\alpha)=r(\beta)\in \{ v_1,\cdots, v_{n_0}\} \Big\},\\
\omega_\Sigma  (n_0,n_1)&:=\Big\{ \sum  s_{\alpha_i}
s_{\beta_i}^*\, \big |\, s_{\alpha_i} s_{\beta_i}^*\in
\omega(n_0,n_1)\Big\}.\end{align*}
  Note that $\omega_\Sigma  (n_0,n_1)$ is not the linear span of $\omega(n_0,n_1)$.
  Then $\{\omega_\Sigma (n_0,n_1)\mid n_0\in \mathbb N,n_1\in \mathbb Z^+\}$ is a
  net of finite subsets in $\mathcal A_E$ which is
  partially ordered by  inclusion. In fact, given two finite sets
  $\omega_\Sigma (n_0,n_1)$, $\omega_\Sigma (m_0,m_1)$ ($n_1\leq
  m_1$), one may write each element $s_\alpha s_\beta^*\in \omega
  (n_0,n_1)$ as $$s_\alpha s_\beta^*=
  s_{\alpha} (\sum_{|\mu|=m_1-n_1} s_\mu s_\mu^*) s_\beta^*=
  \sum  s_{\alpha\mu} (s_{\beta\mu})^*\in \omega_\Sigma
  (m_2,m_1),$$
where $m_2>\max \{n_0, m_0\}$ is an integer large enough so that
$r(\alpha\mu)\in \{v_1, \cdots, v_{m_2}\}$ for any $\alpha\mu$
appearing in the last sum, then clearly $\omega_\Sigma
(n_0,n_1)\cup\omega_\Sigma (m_0,m_1)$ is contained in
$\omega_\Sigma (m_2,m_1).$

   Since the linear span of the set $\cup_{n_0,n_1,n}\Phi_E^n (\omega_\Sigma
   (n_0,n_1))$ is dense in $\mathcal A_E$, by Proposition 3.7, we
   show that for each finite set  $\omega_\Sigma(n_0,n_1)$,
$$ht(\Phi_E,\omega_\Sigma (n_0,n_1))\,\leq\, \max\{h_b({}^t E), h_b(E)\}.$$
If $s_\alpha s_\beta^*\in \omega(n_0,n_1)$,
$r(\alpha)=r(\beta)=v$,  then for $l\leq n-1$,
$$\Phi_E^l(s_\alpha s_\beta^*)=\sum_{|\mu|=l} s_{\mu\alpha}
s_{\mu\beta}^*=\sum_{|\mu|=l}s_{\mu\alpha}(\sum_{\substack{{|\nu|=n-l}\\
{s(\nu)=v}} } s_\nu s_\nu^*) s_{\mu\beta}^*=
\sum_{\substack{{|\mu\alpha\nu|=n+n_1}\\{|\mu|=l}}}
s_{\mu\alpha\nu}
(s_{\mu\beta\nu})^*,$$ because $p_v=\sum_{\substack{{|\nu|=n-l}\\
{s(\nu)=v}} } s_\nu s_\nu^*.$
  Hence one sees  that
$$\cup_{i=0}^{n-1}
\Phi_E^i(\omega_\Sigma(n_0,n_1)) \subseteq\   \Big\{
\sum_{|\mu\alpha\nu|=n+n_1}s_{\mu\alpha\nu}
(s_{\mu\beta\nu})^*\,\big|\, s_\alpha s_\beta^*\in
\omega(n_0,n_1)\Big\}.$$
 Since the set $\{s_{\mu}s_\nu^*\mid \mu,\nu\in
\cup_{i=1}^{n_0} E^{n_1+n}(v_i)\}$ forms a matrix units, it
generates the $C^*$-subalgebra of $\mathcal A_E$ which is
isomorphic to $M_{k_n}$, where $k_n=|\cup_{i=1}^{n_0}
E^{n_1+n}(v_i)|$. Let
 $$\rho_n:{\rm span}\{s_\alpha s_\beta^*\mid \alpha,\beta\in
\cup_{i=1}^{n_0} E^{n_1+n}(v_i)\}\to M_{k_n}$$ be a
$*$-isomorphism with the inverse $\rho^{-1}$. Then by Arveson's
extension theorem $\rho$ extends to a contractive cp map
$\bar{\rho}:\mathcal A_E\to M_{k_n}$, so that we obtain an element
$(\bar{\rho},\rho^{-1},M_{k_n})\in CPA(id, \mathcal A_E)$ such
that $\|\rho^{-1}\circ\bar{\rho}(x)-x\|=0$ if
$$ x\ \in
\cup_{i=0}^{n-1} \Phi_E^i(\omega_\Sigma(n_0,n_1)) \subseteq   \
{\rm span}\{s_\alpha s_\beta^*\mid \alpha,\beta\in
\cup_{i=1}^{n_0} E^{n_1+n}(v_i)\}.$$ Hence
$$rcp\big(\cup_{i=0}^{n-1} \Phi_E^i(\omega_\Sigma(n_0,n_1)),\,\delta\big)
\,\leq \, k_n$$ holds for
any $\delta>0$.
 Thus
$$ht(\Phi_E,\omega_\Sigma (n_0,n_1))\,\leq\, \limsup_{n\to \infty}
\frac{1}{n} \log (k_n).$$
 On the other hand, the irreducibility of $E$ implies that there is
 an $N$ such that  $|E^{n_1+n}(v_i)|\leq |E^{n_1+n+N}(v_1)|$
 for $1\leq i\leq n_0$. Hence $k_n=|\cup_{i=1}^{n_0} E^{n_1+n}(v_i)|\leq
 n_0 |E^{n_1+n+N}(v_1)|.$
Therefore
 $$\limsup_{n\to \infty} \frac{1}{n} \log
 k_n \leq \limsup_{n\to \infty} \frac{1}{n} \log
 |E^{n}(v_1)|,$$ and
the assertion then follows from Proposition 3.5(b).
 \endproof  \medbreak

\vskip 1pc
\begin{ex} \rm Let $E:=E_{\{r_n\},\{l_n\}}$ be a Salama's infinite
irreducible graph (see \cite{Sa}). We assume here that $l_n+1\leq
l_{n+1}$ for each $n$. There are $r_k$ edges from the vertex $k-1$
to $k$, and there is only one path (of length $l_k-l_{k-1}$) from
the vertex $v_k$ to $v_{k-1}$.

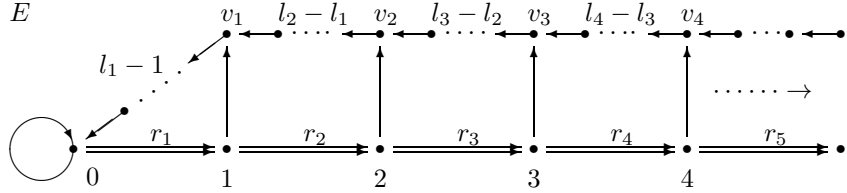
\setlength{\unitlength}{1.7cm}

\hspace*{6cm}
\begin{picture}(8,1.7)

\put(-3.25,0){\circle{0.5}}\put(-2.9,-0.28){$0$}
\put(-3,0){\circle*{0.07}}
 \put(-2.6,0.3){\circle*{0.07}} \put(-2.2,0.55){$\cdot$}
\put(-1.8,0){\circle*{0.07}}\put(-1.85,-0.3){1}
\put(-1.85,1){$v_1$}
\put(-0.6,0){\circle*{0.07}}\put(-0.65,-0.3){$2$}\put(-0.65,1){$v_2$}
\put(0.6,0){\circle*{0.07}}\put(0.55,-0.3){$3$}\put(0.55,1){$v_3$}

\put(1.8,0){\circle*{0.07}}\put(1.75,-0.3){$4$}\put(1.75,1){$v_4$}

\put(3,0){\circle*{0.07}}

 \put(-1.8,0.9){\circle*{0.07}}
\put(-1.4,0.9){\circle*{0.07}}\put(-1.04,0.86){$\cdot$}

\put(-0.6,0.9){\circle*{0.07}}\put(-0.2,0.9){\circle*{0.07}}\put(0.18,0.86){$\cdot$}
\put(0.6,0.9){\circle*{0.07}}
\put(1,0.9){\circle*{0.07}}\put(1.35,0.86){$\cdot$}\put(1.8,0.9){\circle*{0.07}}
\put(2.2,0.9){\circle*{0.07}}\put(2.6,0.9){\circle*{0.07}}\put(3,0.9){\circle*{0.07}}

\put(-1.4,0.9){\vector(-1,0){0.3}} \put(-1.3,0.9){$\dots$}
\put(-0.6,0.9){\vector(-1,0){0.3}}
\put(-0.2,0.9){\vector(-1,0){0.3}} \put(-0.1,0.9){$\dots$}
\put(0.6,0.9){\vector(-1,0){0.3}}

\put(1.0,0.9){\vector(-1,0){0.3}}

\put(1.1,0.9){$\dots$}

\put(1.8,0.9){\vector(-1,0){0.3}}

\put(2.2,0.9){\vector(-1,0){0.3}}

\put(2.3,0.9){$\dots$}

\put(3,0.9){\vector(-1,0){0.3}}

\put(-3.01,0.07){\vector(1,-2){0}}

\put(-2.9,0.02){\vector(1,0){1.0}}\put(-2.9,-0.02){\vector(1,0){1.0}}

\put(-1.7,0.02){\vector(1,0){1.0}}\put(-1.7,-0.02){\vector(1,0){1.0}}

\put(-0.5,0.02){\vector(1,0){1.0}}\put(-0.5,-0.02){\vector(1,0){1.0}}

\put(0.7,0.02){\vector(1,0){1.0}}\put(0.7,-0.02){\vector(1,0){1.0}}

\put(1.9,0.02){\vector(1,0){1.0}}\put(1.9,-0.02){\vector(1,0){1.0}}

\put(-1.8,0.1){\vector(0,1){0.7}}

\put(-0.6,0.1){\vector(0,1){0.7}}

\put(0.6,0.1){\vector(0,1){0.7}}

\put(1.8,0.1){\vector(0,1){0.7}}

\put(-2.6,0.3){\vector(-4,-3){0.3}}

\put(-2.3,0.47){$\cdot$}

\put(-2.4,0.4){$\cdot$}

\put(-2.5,0.33){$\cdot$}

\put(-1.8,0.9){\vector(-4,-3){0.3}}

\put(2,0.4){$\cdots\cdots\rightarrow$}

\put(-2.4,0.06){$r_1$}

\put(-1.2,0.06){$r_2$}\put(0,0.06){$r_3$}\put(1.2,0.06){$r_4$}\put(2.4,0.06){$r_5$}

\put(-3.5,1.0){$E$}

\put(-2.8,0.6){$l_1-1$}

\put(-1.4,1){$l_2-l_1$}

\put(-0.2,1){$l_3-l_2$}

\put(1,1){$l_4-l_3$}
\end{picture}

\vskip 1cm

\noindent Note that for each $n$, $|E^n_r(0^\star)|\leq
|E^n_s(0^\star)|$, which then implies by Proposition 3.2 that
$$h_b({}^tE)\leq h_b(E).$$ Thus
from Theorem 3.9, we have $$ht(\Phi_E|_{\mathcal A_E})\leq
h_b(E).$$
   In particular, if  $E_p:=E_{p,p}$ $(p>1)$ is an
   irreducible infinite graph of Salama satisfying
  $h_l(E_p)=h_b(E_p)=\log p,$  we have
 $$ht(\Phi_{E_p}|_{\mathcal
A_{E_p}})=\log p.$$
\end{ex}

\end{document}